\documentclass{amsart}

\usepackage{graphicx}        
\usepackage{multicol}        

\usepackage{subcaption}
\usepackage{orcidlink}
\usepackage{hyperref}

\newcommand{\pidiv}{\pi^{\text{div}}}

\newcommand{\NSE}{Navier-Stokes equation}

\makeindex             

\title[Gradient-robustness in optimization]{Gradient-robustness in optimization subject to stationary Navier-Stokes equations}

\author[Constanze Neutsch]{Constanze Neutsch\,\orcidlink{0009-0009-5781-8288}}
\address{Department of Mathematics, University of Hamburg, Hamburg, Germany}
\email{constanze.neutsch@uni-hamburg.de}

\author[Winnifried Wollner]{Winnifried Wollner\,\orcidlink{0000-0002-6571-8043}}
\address{Department of Mathematics, University of Hamburg, Hamburg, Germany}
\email{winnifried.wollner@uni-hamburg.de}


\begin{document}

%
%
\begin{abstract}
  In this article, we discuss gradient robust discretizations for the
simulation of the non-linear incompressible Navier-Stokes problem and the optimal
control of such flow. We consider several formulations of the flow
problem that are equivalent for the continuous non-linear forward problem
and compare their gradient robust discretization. We will then discuss
the influence of the chosen formulation on the adjoint equations needed
for gradient computation in the solution of the optimal control problem.
\end{abstract}
\maketitle


\section{Introduction}
\label{sec:1}
For an exterior forcing $\mathbf{f}$, a pair $(\mathbf{u},p) \in (\mathbf{V} \times Q)$,
consisting of a velocity field $\mathbf{u}$ and a pressure $p$  
is said to be a solution to the incompressible stationary \NSE{}, in convective form,
if it satisfies:
\begin{equation} \label{eqn:contiNS}
\begin{aligned}
    \nu (\nabla \mathbf{u}, \nabla \mathbf{\varphi})
    + ((\mathbf{u} \cdot \nabla) \mathbf{u}, \varphi)
    + (\nabla \cdot \varphi, p) &= (\mathbf{f}, \varphi)
    &&\forall \varphi \in \mathbf{V},\\
    (\nabla \cdot \mathbf{u}, \psi) &= 0  &&\forall \psi \in Q
  \end{aligned}                  
\end{equation}            
for function spaces $\mathbf{V}:=\mathbf{H}_0^1(\Omega, \mathbb{R}^d)$ and $Q:=L^2_0(\Omega)$.
Here, $\nu$ denotes the kinematic viscosity, which is inversely related to the Reynolds
number $\operatorname{Re}$ in the non-dimensionalized form of the equation.

The solution of this equation satisfy an important invariance property under gradient `forces'.
More precisely, utilizing the Helmholtz decomposition, any vector field
$\mathbf{f} \in L^2(\Omega)^d$ can be decomposed as
\begin{equation}
    \mathbf{f} = \mathbf{w} + \nabla \varphi \label{eqn:HelmholtzDecomp}
\end{equation}
into an irrotational part $\nabla \varphi$ and a divergence-free part $\mathbf{w}$, i.e., $\nabla \cdot \mathbf{w} = 0$.
Here, the irrotational and the divergence-free part are orthogonal
with respect to the $L_2$-scalar product, i.e.,
\[
(\mathbf{w}, \nabla \varphi)_{L^2} = 0.
\]
Clearly, testing~\eqref{eqn:contiNS} with testfunctions in
$\mathbf{V^0} = \{\mathbf{v} \in \mathbf{V} : \nabla \cdot \mathbf{v} = 0\}$, the solution
$\mathbf{u}$ depends on the divergence-free part $\mathbf{w}$ only, while the irrotational
gradient component needs to be balanced by the pressure. This gives rise to invariance
of the velocity with respect to irrotational forces: if for some arbitrary
forcing $\mathbf{f}$ the pair $(\mathbf{u},p)$ is the  corresponding
solution to~\eqref{eqn:contiNS}, then an irrotational perturbation of the forcing
$\mathbf{f}$ to ($\mathbf{f} + \nabla \varphi$)
evokes a change in the solution pair to be $(\mathbf{u},p+\varphi)$,
see also~\cite{johnDivergenceConstraintMixed2017}.
Pushing the non-linearity to the right hand side shows that similarly the velocity $\mathbf{u}$
is not influenced by gradients in the non-linearity $(\mathbf{u}\cdot\nabla)\mathbf{u}$.

However, applying standard mixed finite elements, this property breaks apart during the
discretization as in general discretely divergence-free functions are not divergence-free.
Therefore, discretely divergence-free functions are not necessarily orthogonal to gradients
in $\mathbf{f}$ or the discrete non-linearity.
Consequently, the discrete velocity is no longer invariant under gradient changes
in the discrete setting~\cite{linkeRoleHelmholtzDecomposition2014}. 
It has been observed that as a consequence spurious spikes and non-physical oscillations
can occur in the computed discrete velocity~\cite{linkeDivergencefreeVelocityReconstruction2012}.
This unfavorable coupling effect is reflected in the error estimate for standard mixed
finite elements in the sense that the velocity error depends on the approximability of both velocity and pressure. Consequently a hard to approximate pressure influences the velocity
error negatively, see, e.g.,~\cite{boffiMixedFiniteElement2013,giraultFiniteElementMethods1986,brezziMixedHybridFinite1991}.

Non-robustness in mixed finite element methods is a well-known fact in
literature~\cite{galvinStabilizingPoorMass2012}.
Substantial research interest has sparked in this area, resulting in several approaches
to gain control over this issue. 
As mentioned, the underlying issue arises from the weakening of the incompressibility
constraint, in the sense that it is enforced only at the discrete level. 
The weakening of the continuity constraint is a price often payed in order to fulfill the
discrete inf-sup condition and have simple to construct spaces~\cite{johnFiniteElementMethods2016}.
The author of~\cite{linkeDivergencefreeVelocityReconstruction2012} proposes the
term \textit{poor momentum balance}, preferred over the more
established \textit{poor mass conservation} (occurring in the divergence-free constraint)
in order to emphasize the incorrect balance between pressure gradients and velocity in the
discrete momentum equation (compared to the continuous case).
Some of the various approaches to circumvent the issue such as, e.g., grad-div
stabilization, e.g.,~\cite{olshanskiiGradDivStabilization2009,olshanskiiGraddivStablilizationStokes2003},
de Rham Complex, $\mathbf{H}$-div conforming spaces, e.g.,~\cite{HdivLederer} and
the use of a interpolation operator are discussed in~\cite{johnDivergenceConstraintMixed2017}. 
Within this paper we will follow the last-mentioned approach of an interpolation operator. 
It was first introduced in~\cite{linkeDivergencefreeVelocityReconstruction2012} and
gave rise to a number of publications building on this idea,
see, e.g.,~\cite{linkeRoleHelmholtzDecomposition2014,linkeQuasioptimalityPressurerobustNonconforming2018,linkeRobustArbitraryOrder2016} also for problems in 
elasticity~\cite{basavaGradientRobustMixed2023,fuLockingFreeGradient2020}. 
The above mentioned publications show promising convergence results, mostly for the
Stokes problem. However, the problem is obviously relevant for
the \NSE{} as well, see, e.g.,~\cite{linkeVelocityErrorsDue2016,johnDivergenceConstraintMixed2017,linkeRoleHelmholtzDecomposition2014}.
As discussed above, in this context, difficulties need not only arise due to the external
force $\mathbf{f}$, but may also be triggered by large gradients stemming from the non-linearity.

The applicability of the interpolation operator approach to \NSE{} is considered
as helpful to also cure poor momentum balance occurring in the non-linear terms of
the \NSE{} and further in the time-derivative in the case of time-dependent
flows~\cite{ahmedPressurerobustDiscretizationOseens2020}. 
The technique has been applied to different Finite Elements (FE) such as Taylor-Hood and
Mini elements in~\cite{ledererDivergencefreeReconstructionOperators2016}
(for Stokes and Navier Stokes), Crouzeix-Raviart~\cite{linkeOptimalL2Velocity2017}, 
and BDM~\cite{brezziTwoFamiliesMixed1985} elements in~\cite{brenneckeOptimalPressureIndependent$L^2$2015}.

Different formulations of the non-linearity in~\eqref{eqn:contiNS} and their error
estimates have been
studied in~\cite{reyesExamplesIdentitiesInequalities2023,linkePressurerobustnessDiscreteHelmholtz2016,charnyiConservationLawsNavier2017}
where the latter focuses on the conservation properties of the different formulations.
The nonlinearities will be detailed in Section~\ref{sec:2} later and
are all equivalent in the continuous case and for exactly divergence-free testfunctions.

However, as with many other ODEs and PDEs, \NSE{} can not just be seen as a stand-alone
equation but can of course be relevant as a constraint in optimization problems, see
~\cite{HinzeHabi,troltzschSecondorderSufficientOptimality2006, vexlerOptimalControlNavierStokes2025, chrysafinosSymmetricErrorEstimates2015, 
casasVelocityTrackingProblem2026, ulbrichConstrainedoptimalControlNaviera}.
There exist promising results of optimal control with (Stokes) and Navier-Stokes
equation~\cite{gunzburgerAnalysisFiniteElement, gunzburgerAnalysisFiniteElement1991} 
\cite{delosreyesOptimalControlStationary2007}(control and state constrained).
Differing in the chosen setting regarding control and or state constraints and
their specific details.
An entirely state constrained case by stationary \NSE{} be found in~\cite{delosreyesOptimalControlStationary2007}.

The concept of gradient robust discretizations for control problems subject to Stokes flow
has only recently been considered in~\cite{merdonPressureRobustnessContextOptimal2023}. 
In particular, it was highlighted that the possibly occurring irrotational forces in the
cost functional can similarly disturb the velocity calculation via the adjoint equation.
Hence a gradient robust discretization of an optimal control problem needs to consider
not only a gradient robust discretization of the state equation but also the resulting
optimality system which introduces additional sources of non-robustness.
Recent work of~\cite{liPressurerobustnessStokesDarcyOptimal2025} investigated the same
issue of gradient-robustness in the context of optimal control, 
considering a Stokes-Darcy flow as the underlying constraint.

These promising results provide the incentive for this work: extending the gradient-robust
methodology for PDE-constrained optimization subject to 
stationary incompressible \NSE{}.

The focus of this work is given by the following optimization problem,
with the control space $\mathbf{Q}= \mathbf{L}^2(\Omega, \mathbb{R}^n)$:
\begin{equation}\label{eqn:OCP}
\begin{aligned}
\min_{(\mathbf{q},\mathbf{u},p) \in \mathbf{Q}\times \mathbf{V} \times Q } \frac{1}{2} &\lvert \lvert \mathbf{u}- \mathbf{u}^d \rvert \rvert^2_{\mathbf{L}^2(\Omega)} + 
\frac{1}{2} \lvert \lvert \mathbf{q} \rvert \rvert^2_{\mathbf{L}^2(\Omega)} \\
\text{s.t.}\quad \text{ Navier-Stokes equation with } &\text{additional control } \mathbf{q} \text{ on the right hand side.}
\end{aligned}
\end{equation}
The structure of this work is as follows: 
Section~\ref{sec:2} describes the different formulations of the \NSE{} to be used as a
constraint in~\eqref{eqn:OCP}. 
The subsequent Section~\ref{sec:3} deals with the discretization and introduces the
optimization problem constrained by Navier-Stokes equation. 
The numerical results comparing the advantageous robust and non-robust approaches are
presented and analyzed in Section~\ref{sec:4}. 

\section{Navier-Stokes equation}
\label{sec:2}
Having introduced the basic idea of the gradient-robustness, the \NSE{} with its different
formulations of the non-linearity is described in more detail below.
Therefore, the Helmholtz projection $\mathbb{P}\colon L^2(\Omega)^d \rightarrow L^2(\Omega)^d$ is introduced as an addition to the
Helmholtz decomposition~\eqref{eqn:HelmholtzDecomp}, see, e.g.,~\cite{sohrNavierStokesEquations2001}:
\begin{equation}
    \mathbb{P}:L^2(\Omega)^d \rightarrow L^2(\Omega)^d, \mathbf{f} \mapsto \mathbf{w} \label{eqn:HelmholtzProj}\\
\end{equation}
where $\mathbf{f} = \nabla \varphi + \mathbf{w}$ and $\nabla \cdot \mathbf{w}=0$.
The projection operator is essential for the understanding of the robustness, which will
be discussed in further detail shortly.

Considering the convective form $(\mathbf{u} \cdot \nabla)\mathbf{u}$ of the non-linearity
in~\eqref{eqn:contiNS} is not the only possibility. Indeed there are some equivalent
formulations that we would like to consider, in particular for the potential differences
in view of the adjoint equation in the optimality conditions of~\eqref{eqn:OCP}.
For ease of reading we introduce the trilinear form
$c(\cdot, \cdot, \cdot) \colon \mathbf{V}\times \mathbf{V}\times \mathbf{V} \rightarrow \mathbb{R}$ for the representation of the non-linearity in weak form.
For $\mathbf{u}, \mathbf{v}, \mathbf{w}\in \mathbf{V}$, $c(\cdot, \cdot, \cdot)$ encompasses
the following equivalent variants:
\begin{align}
    c_{\text{conv}}(\mathbf{u},\mathbf{v},\mathbf{w}) &= \int ((\mathbf{u} \cdot \nabla)\mathbf{v}) \cdot \mathbf{w} \,\mathrm{d}x, \label{eqn:c_conv} \\ 
    c_{\text{div}}(\mathbf{u},\mathbf{v},\mathbf{w}) &= \int ((\mathbf{u} \cdot \nabla)\mathbf{v}) \cdot \mathbf{w}  + \frac{1}{2} ((\nabla \cdot \mathbf{u})\mathbf{v}) \cdot \mathbf{w}  \,\mathrm{d}x, \label{eqn:c_div}\\ 
    c_{\text{rot}}(\mathbf{u},\mathbf{v},\mathbf{w}) &= \int ((\nabla \times \mathbf{u})\times \mathbf{v})\cdot \mathbf{w} \,\mathrm{d}x,  \label{eqn:c_rot} \\
    c_{\text{skew}}(\mathbf{u},\mathbf{v},\mathbf{w}) &= \int \frac{1}{2} \Big( ((\mathbf{u} \cdot \nabla)\mathbf{v})\cdot \mathbf{w}- ((\mathbf{u} \cdot \nabla)\mathbf{w})\cdot \mathbf{v}  \Big) \,\mathrm{d}x\label{eqn:c_skew}\\ 
     &= \frac{1}{2} (c_{\text{conv}}(\mathbf{u},\mathbf{v},\mathbf{w}) - c_{\text{conv}}(\mathbf{u},\mathbf{w},\mathbf{v})). \nonumber
\end{align}
We would like to begin with a few comments regarding the different formulations of the non-linearity.
\subsubsection*{The convective form~\eqref{eqn:c_conv}}
We would like to highlight that for
any $u = \nabla \chi$ being a gradient of a potential flow 
it holds that $((\mathbf{u} \cdot \nabla) \mathbf{u}) = \frac{1}{2}\nabla (\lvert \mathbf{u} \rvert)^2 $ is also a gradient and consequently
\[
\int_{\Omega} ((\mathbf{u} \cdot \nabla) \mathbf{u}) \mathbf{w} \,\mathrm{d}x = 0 \quad \forall \mathbf{w} \in \mathbf{V}^0
\] 
due to the orthogonality of gradient fields and divergence-free functions of
$\mathbf{V}^0$ with respect to the $L^2$ scalar product. 
For a proof the reader is referred to~\cite{linkePressurerobustnessDiscreteHelmholtz2016},
the following details are likewise taken from that source.
In a more general setting (of not necessarily potential flows $\mathbf{u}$) the convective
term can be challenging in two ways:
Firstly, the more obvious one in the context of gradient-robustness arises from the poor
treatment and breakdown of $L^2$-orthogonality between discretely divergence-free test functions 
and gradient forces.
Similarly to the issue on the right hand side of the equation
($\int \mathbf{f} \cdot \mathbf{v}$), a large gradient part
of $((\mathbf{u} \cdot \nabla) \mathbf{u})$ can induce errors 
when multiplied with a discretely divergence-free function
$\mathbf{w}_h$ as $((\mathbf{u}_h \cdot \nabla) \mathbf{u}_h) \mathbf{w}_h$.
The gradient part can be quantified by $ (\mathbf{u} \cdot \nabla)\mathbf{u} - \mathbb{P}(\mathbf{u} \cdot \nabla)\mathbf{u}$ (recalling the Helmholtz projection~\eqref{eqn:HelmholtzProj}).

A second source of error gets excited for high Reynolds numbers $\operatorname{Re}$ when the
flow behavior switches to a potentially more turbulent state.
In this case the divergence-free part $\mathbb{P}((\mathbf{u} \cdot \nabla)\mathbf{u})$ is
large. Such an advection dominance can be challenging as it may introduce oscillations to
the solution\label{sec:advection-dominance}
\cite{linkePressurerobustnessDiscreteHelmholtz2016}.
To obtain skew-symmetry of the convective form one would need $\mathbf{u}$ to be weakly
divergence-free and $\mathbf{u} \cdot \mathbf{n} = 0 \text{ on } \partial \Omega$,
then it holds $c_{\text{conv}}(\mathbf{u},\mathbf{v},\mathbf{v}) = 0$.
\subsubsection*{The divergence form~\eqref{eqn:c_div}}
The form~\eqref{eqn:c_div} has one additional term compared to the convective term with
a factor of $\frac{1}{2}$ which is motivated by the need for
skew-symmetry. In particular, $c_{\text{div}}(u,\mathbf{v},\mathbf{v}) = 0$ even if $\mathbf{u}$
is not weakly
divergence-free~\cite{johnFiniteElementMethods2016}. 
\subsubsection*{The rotational form~\eqref{eqn:c_rot}}
We would like to point out that 
$c_{\text{conv}}(\mathbf{u}, \mathbf{u}, \cdot) = c_{\text{rot}}(\mathbf{u}, \mathbf{u}, \cdot) + \int_\Omega \frac{1}{2}\lvert \mathbf{u} \vert^2 \cdot\,\mathrm{d}x$, where the latter can be 
merged with the pressure to the Bernoulli pressure $P = p + \frac{1}{2}|\mathbf{u}|^2$.
Consequently, $p$ is replaced by $P$ in~\eqref{eqn:contiNS}, such that:
$ \nu \Delta \mathbf{u} + (\nabla \times \mathbf{u}) \times \mathbf{u}  + \nabla P = f $.
The rotational form is skew-symmetric in the last two arguments and thus $c_{\text{rot}}(\mathbf{u},\mathbf{v},\mathbf{v}) = 0$ for $\mathbf{u},\mathbf{v} \in \mathbf{V}$.
For further remarks the reader is referred to~\cite{laytonAccuracyRotationForm2009}.
\subsubsection*{The skew-symmetric form~\eqref{eqn:c_skew}}
A skew-symmetric modification of the convective form~\eqref{eqn:c_conv} is~\eqref{eqn:c_skew}.
Analogous to the rotational form this implies 
$c_{\text{skew}}(\mathbf{u},\mathbf{v},\mathbf{v}) = 0$ for $\mathbf{u},\mathbf{v} \in \mathbf{V}$. However, adapting this form for robustness is not possible, due to
$\pidiv(\mathbf{v}_h) \notin\mathbf{V}$, see, e.g.,~\cite{linkePressurerobustnessDiscreteHelmholtz2016}.
Therefore, we will not compare this form in the following.

For details regarding stability, including the fact of skew-symmetry, the topic of existence
of solutions and the small data assumption the 
reader is referred to standard literature such as~\cite{johnFiniteElementMethods2016}. 

In the following sections the three first-named formulations~\eqref{eqn:c_conv},~\eqref{eqn:c_div},~\eqref{eqn:c_rot} will be employed.
Thus the weak formulation can be written as  
\begin{equation}\label{eqn:genNS}
\begin{aligned}
    -\nu (\nabla \mathbf{u}, \nabla \mathbf{\varphi}) + c(\mathbf{u},\mathbf{u},\varphi) + (p, \nabla \cdot \varphi)&= (\mathbf{f},\varphi) && \forall \varphi \in \mathbf{V}, \\
    (\nabla \cdot \mathbf{u},\psi) &= 0  && \forall \psi \in Q.
\end{aligned} 
\end{equation}
Here, $c(\mathbf{u},\mathbf{u},\varphi)$ is a placeholder for one of the possible formulations above.
Regarding existence and uniqueness of solutions the reader is referred to~\cite[Chapter IV, Section §2, Theorem 2.1, 2.2]{giraultFiniteElementMethods1986}.

In order to wrap up the Navier-Stokes section with a connection to optimization the optimality system in the continuous case
shall be noted here.
Utilizing standard calculus for any local minimizer
$(\mathbf{q},\mathbf{u},p)$ of~\eqref{eqn:OCP} there is an adjoint
state $(\mathbf{z},s) \in \mathbf{V} \times Q$ solving the 
system: 
\begin{equation}\label{eqn:optisys_conti}
\begin{aligned}
    \nu (\nabla \mathbf{u}, \nabla \varphi) + c(\mathbf{u}, \mathbf{u},\varphi) + (p, \nabla \cdot \varphi) &= (\mathbf{q} +\mathbf{f}, \varphi)  &&\forall \varphi \in \mathbf{V}, \\
    (\nabla \cdot \mathbf{u}, \psi) &= 0 &&\forall \psi \in Q,\\
    \nu (\nabla \varphi, \nabla \mathbf{z}) + \tilde{c}(\mathbf{u}, \varphi, \mathbf{z}) + (\mathbf{z}, \nabla \cdot \varphi) &= (\pidiv(\mathbf{u})-\mathbf{u}^d,\varphi) && \forall \phi\in \mathbf{V}, \\
    (\nabla \cdot \psi, s) &= 0 &&  \forall \psi \in Q,\\
     (\mathbf{q}, \varphi)+(\mathbf{z},\varphi) &=0  && \forall
     \varphi\in \mathbf{Q}
   \end{aligned}
 \end{equation}
where
\[
  \tilde{c}(\mathbf{u}, \varphi, \mathbf{z})=c(\mathbf{u}, \varphi, \mathbf{z})+c(\varphi,\mathbf{u}, \mathbf{z})
\]
is the derivative, with respect to the $\mathbf{u}$ variable in
direction $\varphi$ in the point $(\mathbf{u},\mathbf{z}$), of the
mapping $(\mathbf{u},\mathbf{z}) \mapsto c(\mathbf{u},\mathbf{u},\mathbf{z})$. See also \eqref{c_tilde_discrete} and \eqref{eq:c_tilde_formen} for the concrete form in the discretized robust setup.

\section{Discretization}
\label{sec:3}
For the discretization finite dimensional spaces $\mathbf{V}_h \subset \mathbf{V}, Q_h \subset Q$ 
are considered where we assume that the pair is chosen such that the inf-sup condition
is satisfied.
This pair gives rise to a set of discretely divergence-free functions
\[
        \mathbf{V}_h^{0} :=\{ \varphi_h \in \mathbf{V}_h : (\nabla \cdot \varphi_h, \psi_h) = 0 \quad \forall \psi_h \in Q_h\}
\]
where in general $\mathbf{V}_h^0 \not\subset \mathbf{V}^0$.
Along with these spaces, the discrete \NSE{} can be written as
\begin{equation}\label{eqn:discNSE}
\begin{aligned}
    -\nu (\nabla \mathbf{u}_h, \nabla \varphi_h) + c(\mathbf{u}_h,\mathbf{u}_h,\varphi_h)
    + (p_h,\nabla \cdot \varphi_h)&= (\mathbf{f},\varphi_h) && \forall  \varphi_h \in \mathbf{V}_h, \\
    (\nabla \cdot \mathbf{u}_h,\psi_h) &= 0 && \forall \psi_h \in Q_h.
    \end{aligned}
\end{equation}      
where again the trilinear-form can be any of the forms~\eqref{eqn:c_conv},~\eqref{eqn:c_div},~\eqref{eqn:c_rot}.
The theory for existence stability of solutions is covered in~\cite{linkePressurerobustnessDiscreteHelmholtz2016}.

As already mentioned, the above discrete Navier-Stokes equation is not robust with respect to
gradient forces. We will shortly outline how to obtain a gradient robust discretization, following~\cite{linkeRoleHelmholtzDecomposition2014}.
To avoid lengthy abstract conditions on the involved reconstruction operator $\pidiv$ mapping
$\mathbf{V}_h^0$ into $\mathbf{V^0}$ we fix our discretization. To this end, we assume our
domain is subdivided into quadrilateral (or hexahedral) elements. Then we discretize
the space of velocities $\mathbf{V}_h$ by conforming biquadratic ($\mathcal{Q}_2$) finite
elements while $Q_h$ is discretized by discontinuous linear finite elements ($\mathcal{DGP}_1$).
Then~\cite[Sec.~4.2.1]{linkeRobustArbitraryOrder2016} shows that the operator 
$\pidiv$ given by the canonical interpolation into the Brezzi-Douglas-Marini space of
order two ($\mathcal{BDM}_2$) can be used.
The reader is referred to~\cite{linkeRobustArbitraryOrder2016}
and~\cite{johnDivergenceConstraintMixed2017} for details regarding the operator in the
sense of optimality and error estimates.

As already mentioned in Section~\ref{sec:2} the Helmholtz
decomposition~\eqref{eqn:HelmholtzDecomp} and Helmholtz
projection~\eqref{eqn:HelmholtzProj} are closely 
linked to the concept of the interpolation operator $\pidiv$. 
This shall be shown by the example of the treatment of the forcing term $\mathbf{f}$ on
the right hand side. 
By the construction of the operator and the orthogonality of the divergence-free vector
fields $\mathbf{v} \in \mathbf{V}^0$ and the irrotational fields $\nabla \varphi$
one is allowed to insert $\mathbb{P}$ into the following integrals for the continuous case:
\[
\int_{\Omega} \mathbb{P}(\mathbf{f}) \cdot \mathbf{v} \,\mathrm{d}x = \int_{\Omega} \mathbf{f} \cdot \mathbf{v} \,\mathrm{d}x.
\]
However, this equality does not hold for discretely divergence-free functions
$\mathbf{v}_h \in \mathbf{V}_h^0$, instead, one needs to incorporate the interpolation
operator $\pidiv$
\[
\int_{\Omega} \mathbb{P}(\mathbf{f}) \cdot \mathbf{v}_h \,\mathrm{d}x  = \int_{\Omega} \mathbf{f} \cdot \pidiv(\mathbf{v}_h) \,\mathrm{d}x.
\]
By the same argumentation and the details mentioned in Section~\ref{sec:advection-dominance}
regarding the convective term one obtains:
\[
\int_{\Omega} \mathbb{P}((\mathbf{u}_h \cdot \nabla)\mathbf{u}_h) \cdot\mathbf{v}_h \,\mathrm{d}x \neq \int_{\Omega} ((\mathbf{u}_h \cdot \nabla)\mathbf{u}_h) \cdot \mathbf{v}_h \,\mathrm{d}x,
\] 
but 
\[
\int_{\Omega} \mathbb{P}((\mathbf{u}_h \cdot \nabla)\mathbf{u}_h) \cdot \mathbf{v}_h \,\mathrm{d}x = \int_{\Omega} ((\mathbf{u}_h \cdot \nabla)\mathbf{u}_h) \cdot \pidiv(\mathbf{v}_h) \,\mathrm{d}x.
\] 

As one can already infer from the previous equation, the robust modified version of the
convective term is given by:
\begin{equation}\label{eqn:c_conv_h}
c_{\text{conv}}^h(\mathbf{u}_h,\mathbf{w}_h,\varphi_h) = \int_{\Omega} ((\mathbf{u}_h \cdot \nabla)\mathbf{w}_h) \cdot \pidiv(\varphi_h) \,\mathrm{d}x.
\end{equation}
This way the orthogonality between the gradient arising from
$(\mathbf{u} \cdot \nabla)\mathbf{u}$ is (then divergence-free)
testfunction $\pidiv(\mathbf{v}_h)$ is ensured. 
However, the convective term is, in general, not skew-symmetric after discretization,
since $\mathbf{u}_h$ is only discretely divergence-free.
The robust formulation of the divergence form includes the robust formulation of the
convective term as the first summand:
\begin{equation}\label{eqn:c_div_h}
c_{\text{div}}^h(\mathbf{u}_h,\mathbf{w}_h,\varphi_h) =
 \int_{\Omega} (\mathbf{u}_h \cdot \nabla)\mathbf{w}_h \cdot \pidiv(\varphi_h)  + \frac{1}{2} (((\nabla \cdot \mathbf{u})\mathbf{w}) \cdot \pidiv(\varphi_h) ) \,\mathrm{d}x.
\end{equation}
For the rotational form the robust version is given by:
\begin{equation}\label{eqn:c_rot_h}
c_{\text{rot}}^h(\mathbf{u}_h,\mathbf{w}_h,\varphi_h)  = \int_{\Omega} ((\mathbf{u}_h \times \pidiv(\mathbf{w}_h))\times \pidiv(\varphi_h)).
\end{equation}
By the additional application of $\pidiv$ on the second argument, compared to just the third, skew-symmetry is obtained~\cite{linkeRoleHelmholtzDecomposition2014}.
As the arguments to which the operator $\pidiv$ is applied we abbreviate the discrete
non-linearity by $c^h$ and obtain the gradient robust discrete formulation:
\begin{equation}\label{eqn:discNSErob}
\begin{aligned}
    -\nu (\nabla \mathbf{u}_h, \nabla \varphi_h) + c^h(\mathbf{u}_h,\mathbf{u}_h,\varphi_h)
    + (p_h,\nabla \cdot \varphi_h)&= (\mathbf{f},\pidiv(\varphi_h)) && \forall  \varphi_h \in \mathbf{V}_h \\
    (\nabla \cdot \mathbf{u}_h,\psi_h) &= 0 && \forall \psi_h \in Q_h.
\end{aligned}
\end{equation}

In the context of optimal control with the (Navier-)Stokes equation as constraint,
robustness must be ensured not only in the 
state equation but also in the cost functional to properly treat
occurring gradient forces in the adjoint equation.
This is referred to as \textit{fully pressure-robust scheme}
in~\cite{merdonPressureRobustnessContextOptimal2023},
showing that it comes along with the advantage of having a robust adjoint
(for the Stokes equation),
compared to the \textit{partially pressure-robust scheme}; where only the state equation is
modified.
The robust version of the tracking type cost functional with
$\mathbf{u}^d \in L^2(\Omega)$ being the desired solution is then given by:
 \begin{equation}\label{eqn:discOCProb}
  \min_{(\mathbf{q}_h, \mathbf{u}_h, p_h) \in \mathbf{Q}_h, \mathbf{V}_h, Q_h} \frac{1}{2} \lvert \lvert \pidiv(\mathbf{u}_h) - \mathbf{u}^d \rvert \rvert^2_{L^2(\Omega)} + \frac{1}{2} \lvert \lvert \mathbf{q}_h \rvert \rvert^2_{L^2(\Omega)} \qquad
  \text{s.t.~\eqref{eqn:discNSErob}}
\end{equation}
where again $\mathbf{f} = \mathbf{f}+\mathbf{q}_h$. The discrete control space $\mathbf{Q}_h$
is chosen as $\mathbf{Q}_h = \mathbf{V}_h$, which coincides with the variational
discretization of~\cite{hinzeVariationalDiscretizationConcept2005}. \label{sec:vari_discri}

For a local minimizer $(\mathbf{q}_h,\mathbf{u}_h,p_h)$ of~\eqref{eqn:discOCProb}, standard calculus yields the necessary optimality conditions, with an adjoint state
$(\mathbf{z}_h,s_h) \in \mathbf{V}_h \times Q_h$:
\begin{equation}\label{eqn:optisys_discrete} 
\begin{aligned}
    \nu (\nabla \mathbf{u}_h, \nabla \varphi_h) + c^h(\mathbf{u}_h, \mathbf{u}_h,\varphi_h) + (p, \nabla \cdot \varphi_h) &= (\mathbf{q}_h+\mathbf{f}, \varphi_h)  &&\forall \varphi_h \in \mathbf{V}_h, \\
    (\nabla \cdot \mathbf{u}_h, \psi_h) &= 0 &&\forall \psi_h \in Q_h,\\
    \nu (\nabla \varphi_h, \nabla \mathbf{z}_h) + \tilde{c}^h(\mathbf{u}_h, \varphi_h, \mathbf{z}_h) + (\mathbf{z}_h, \nabla \cdot \varphi_h) &= (\pidiv(\mathbf{u}_h)-\mathbf{u}^d,\varphi_h) && \forall \phi_h\in \mathbf{V}_h, \\
    (\nabla \cdot \psi_h, s) &= 0 &&  \forall \psi_h \in Q_h,\\
     (\mathbf{q}_h, \varphi_h)+(\mathbf{z}_h,\varphi_h) &=0  && \forall \varphi_h\in \mathbf{Q}_h  
   \end{aligned}
 \end{equation}
where
\begin{equation}
    \tilde{c}^h(\mathbf{u}_h,\varphi_h,  \mathbf{z}_h) := c^h(\varphi_h, \mathbf{u}_h, \mathbf{z}_h) + c^h(\mathbf{u}_h, \varphi_h, \mathbf{z}_h). \label{c_tilde_discrete}
\end{equation}
The corresponding expressions for $\tilde{c}^h$ for the three formulations under consideration are given by:
\begin{equation}\label{eq:c_tilde_formen}
\begin{aligned}
    \tilde{c}^h_{\text{conv}}(\mathbf{u}_h,\varphi_h,  \mathbf{z}_h) &= ((\mathbf{u}_h \cdot \nabla)\cdot \varphi_h, \pidiv(\mathbf{z}_h)) + ((\varphi_h \cdot \nabla)\cdot \mathbf{u_h}, \pidiv(\mathbf{z}_h),) \\
    \tilde{c}^h_{\text{div}}(\mathbf{u}_h,\varphi_h,  \mathbf{z}_h) &= ((\mathbf{u}_h \cdot \nabla)\cdot \varphi_h, \pidiv(\mathbf{z}_h)) + \frac{1}{2}((\nabla \cdot \mathbf{u}_h)\cdot \varphi_h, \pidiv(\mathbf{z}_h)) \\ 
                                                                    &\;\;\;\; +((\varphi_h \cdot \nabla)\cdot \mathbf{u_h}, \pidiv(\mathbf{z}_h)) + \frac{1}{2}((\nabla \cdot \varphi_h)\cdot \mathbf{u}_h, \pidiv(\mathbf{z}_h)),\\
    \tilde{c}^h_{\text{rot}}(\mathbf{u}_h,\varphi_h,  \mathbf{z}_h) &=
    ((\nabla \times \mathbf{u}_h)\times \pidiv(\varphi),
    \pidiv(\mathbf{z_h})) \\
    &\;\;\;\;+ ((\nabla \times \varphi_h)\times \pidiv(\mathbf{u_h}), \pidiv(\mathbf{z})).
\end{aligned}
\end{equation}

\section{Numerical studies}
\label{sec:4}
The numerical examples were obtained with the finite element library deali.ii \cite{2025:arndt.bangerth.ea:deal} 
along with the finite element/optimization library DOpElib \cite{DOpElib}.
The examples serve the purpose of examplarily identifying the applicability of the
interpolation operator $\pidiv$ in the context of optimal 
control under the constraint of an incompressible \NSE{} with different formulations of the
non-linearity (convective~\eqref{eqn:c_conv}, divergence~\eqref{eqn:c_div} and rotational~\eqref{eqn:c_rot}), as aforementioned.
The error of interest is given by the $L_2$ norm of the gradient of the deviation between the
calculated discretized optimal state $\mathbf{u}_h$
from~\eqref{eqn:optisys_discrete} and the optimal state
$\mathbf{u}$ given by~\eqref{eqn:optisys_conti}, which
is known analytically, here:
$
\lvert \lvert \nabla(  \mathbf{u} - \mathbf{u}_h) \rvert \rvert_{L^2(\Omega)}$. 
The domain $\Omega$ is given by $\Omega = (-1,1)^2$. 
Further, the example is constructed as follows: Let $\mathbf{u} =
\nabla \chi(x,y) = \nabla (x^3 -3xy^2)$. It is taken
from~\cite[Example~1, 2]{linkePressurerobustnessDiscreteHelmholtz2016}
neglecting time-dependence and by
definition irrotational. Further, the desired state is taken as a
gradient perturbation of the optimal state as 
$\mathbf{u}^d = \mathbf{u} + \nabla \psi$ with $\psi = (-(10(x - 0.5)^3 y^2 + (1 - x)^3 (y - 0.5)^3 + 1/8))$ (taken from \cite{basavaGradientRobustMixed2023}).
As a consequence, the optimal adjoint satisfies $\mathbf{z} = 0$ while
the adjoint pressure $s$ balances the additional gradient in the cost
functional.  Therefore, a second error functional 
$\lvert \lvert \nabla \mathbf{z}_h \rvert \rvert $ is used to measure
the error in the adjoint velocity.
The boundary conditions are given in terms of Dirichlet values of the solution $\mathbf{u}$. 
The forcing is set to $\mathbf{f} = 0$, along with an initialization of the control $\mathbf{q}_h = 0$.

For the calculations a rectangular uniform grid has been used, equipped with
$\mathbf{V}_h = Q_2$ elements for the velocity, $\mathbf{Q}_h = Q_2$ for the control and
$Q_h = \mathcal{DGP}_1$ for the pressure. 
The concept of variational discretization is realized by $\mathbf{Q}_h = \mathbf{V}_h$ (as mentioned in Section \ref{sec:vari_discri}).
The interpolation operator utilizes $\mathcal{BDM}_2$ elements as a projection
space~\cite{brezziTwoFamiliesMixed1985}. 
The viscosity $\nu$ is varied over three orders of magnitude
taking the values $\{$1.0, 0.1, 0.01$\}$.
The problem is solved on a uniform rectangular grid with 722,946 DoFs for the state and 132,098 DoFs for the control.
 
The construction of $\mathbf{u}^d$ and $\mathbf{u}$ makes the computation theoretically susceptible to 
suffer from errors when not treated robustly, as the term $(\mathbf{u}-\mathbf{u}^d)$  of the cost functional 
is a gradient and is present in the adjoint equation on the right hand
side. Thus, for the adjoint $\mathbf{z}_h\neq 0$ is expected when a
non-robust discretization of the adjoint equation is used.
Theoretically, one expects the non-robust schemes to have errors larger than the respective
robust scheme and additionally observe the error scaling with the order of magnitude of the viscosity $\nu$.
In Table~\ref{tab:results_PiImFunktionalQ2UndRobQFormen} the results of the three formulations in the robust and non-robust case
with respect to the gradient error are compared.
\begin{table}[h!]
\caption{Error $\lvert \lvert \nabla (\mathbf{u}^d - \mathbf{u}_h) \rvert \rvert_{L^2(\Omega)}$ of robust and non-robust method, using the \textit{convective form}, \textit{divergence form} 
and \textit{rotational form }for varying viscosities $\nu$.} 
\begin{tabular}{p{0.12\textwidth}|p{0.13\textwidth}p{0.13\textwidth}|
                        p{0.13\textwidth}p{0.13\textwidth}|
                        p{0.13\textwidth}p{0.13\textwidth}}
    \hline\noalign{\smallskip}
& \multicolumn{2}{c|}{\parbox{0.26\textwidth}{\centering Convective}}
& \multicolumn{2}{c|}{\parbox{0.26\textwidth}{\centering Divergence}}
& \multicolumn{2}{c}{\parbox{0.26\textwidth}{\centering Rotational}} \\
\noalign{\smallskip}\hline\noalign{\smallskip}
Viscosity $\nu$ 
&  Robust & Non-Robust & Robust & Non-Robust & Robust & Non-Robust \\
\noalign{\smallskip}\hline\noalign{\smallskip}
1e+0 & 2.307e-13 & 1.066e-6 & 2.284e-13 & 1.066e-6 & 2.334e-13 & 2.353e-13 \\
1e-1 & 2.362e-13 & 1.066e-5 & 2.390e-13 & 1.066e-5 & 2.398e-13 & 2.384e-13 \\
1e-2 & 3.896e-13 & 1.066e-4 &3.855e-13  & 1.066e-4 & 2.459e-13 & 2.429e-13 \\
\noalign{\smallskip}\hline\noalign{\smallskip}
\end{tabular}\\
\label{tab:results_PiImFunktionalQ2UndRobQFormen}
\end{table}
One can deduce that the robust approach is
highly advantageous for the convective and divergence formulation
compared to the respective non-robust schemes. 
Additionally, the order of the error of the non-robust schemes scales with the order of the
viscosity, which is in line with established error estimates for the
discretization of the forward problem, see, e.g.,~\cite{linkeRoleHelmholtzDecomposition2014,linkePressurerobustnessDiscreteHelmholtz2016}.
The robust schemes are almost entirely viscosity-independent, which meets the expectation. 
Regarding the rotational form almost no differences in accuracy can be observed for the robust or non-robust approach. 
This can be attributed to the fact the nonlinearity vanishes in the
forward equation for this particular choice of $\mathbf{u}$ and the
absence of volume forces $\mathbf{f}$.
However, for the adjoint equation both a volume force
$\mathbf{u}-\mathbf{u}^d$ is present and the linearized non-linearity
does not vanish. Consequently differences in the adjoint accuracy
(compared to the correct solution of $\mathbf{z}=0$) can be observed
for the non-robust discretization, 
see Table~\ref{tab:nablaz}.
\begin{table}[h!]
\caption{Error $\lvert \lvert \nabla \mathbf{z}_h \rvert \rvert_{L^2(\Omega)}$ of robust and non-robust method, using the \textit{convective form}, \textit{divergence form} 
and \textit{rotational form }for varying viscosities $\nu$.}
\begin{tabular}{p{0.12\textwidth}|p{0.13\textwidth}p{0.13\textwidth}|
                        p{0.13\textwidth}p{0.13\textwidth}|
                        p{0.13\textwidth}p{0.13\textwidth}}
    \hline\noalign{\smallskip}
& \multicolumn{2}{c|}{\parbox{0.26\textwidth}{\centering Convective}}
& \multicolumn{2}{c|}{\parbox{0.26\textwidth}{\centering Divergence}}
& \multicolumn{2}{c}{\parbox{0.26\textwidth}{\centering Rotational}} \\
\noalign{\smallskip}\hline\noalign{\smallskip}
Viscosity $\nu$ 
&  Robust & Non-Robust & Robust & Non-Robust & Robust & Non-Robust \\
\noalign{\smallskip}\hline\noalign{\smallskip}
1e+0 & 2.600e-8 & 3.502e-6 & 4.081e-8 & 3.502e-6 & 2.602e-8 & 3.502e-6\\
1e-1 & 2.600e-8 & 3.502e-5 & 4.082e-8 & 3.502e-5 & 2.602e-8 & 3.502e-5\\
1e-2 & 2.600e-8 & 3.501e-4 & 4.082e-8 & 3.501e-4 & 2.602e-8 & 3.502e-4\\
\noalign{\smallskip}\hline
\end{tabular}\\
\label{tab:nablaz}
\end{table}

\begin{figure*}[h!]
    \centering
    \begin{subfigure}[t]{0.3\textwidth}
        \includegraphics[scale=0.13, trim=90pt 0 150pt 0, clip]{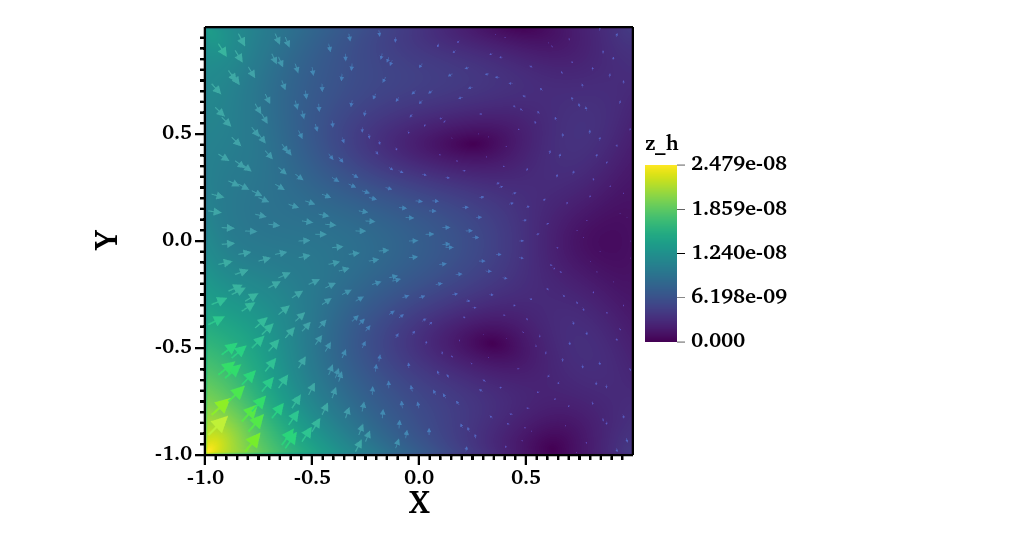}       
         \caption{$\nu = $ 1e+0, non-robust}
       \end{subfigure}%
       \hfill
       \begin{subfigure}[t]{0.3\textwidth}
         \centering
         \includegraphics[scale=0.13, trim=150pt 0 150pt 0, clip]{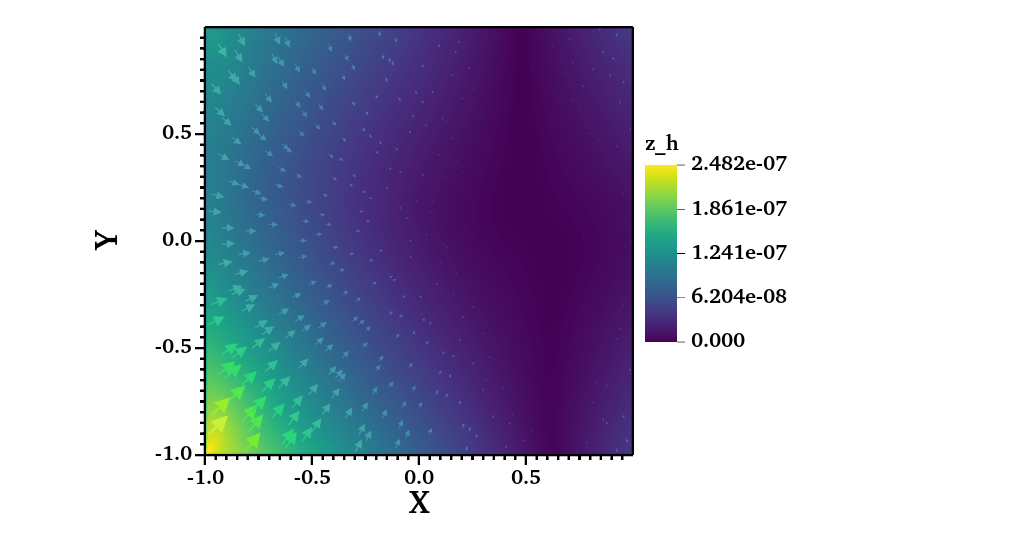}
         \caption{$\nu = $ 1e-1, non-robust}
       \end{subfigure}
       \hfill
       \begin{subfigure}[t]{0.3\textwidth}
         \centering
         \includegraphics[scale=0.13, trim=150pt 0 150pt 0, clip]{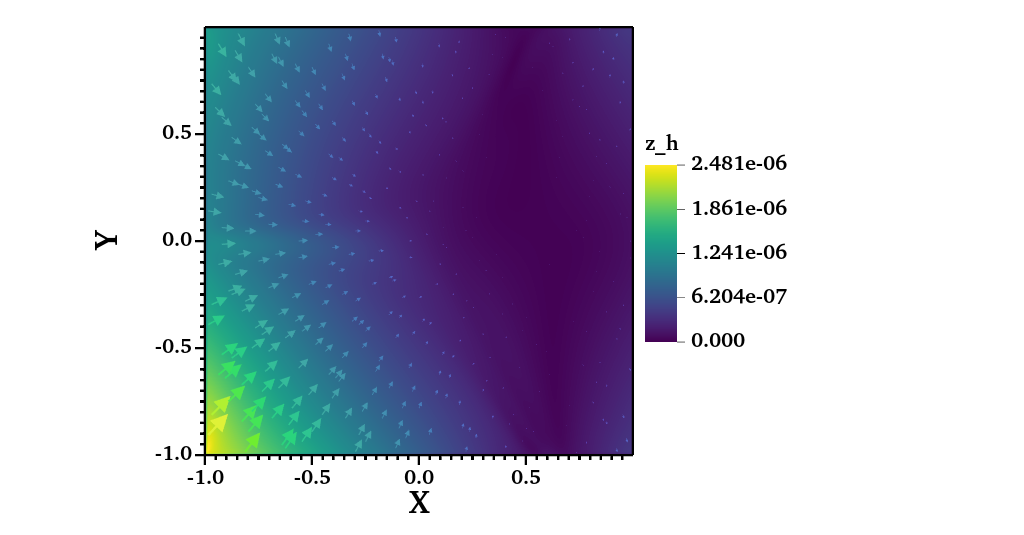}
         \caption{$\nu = $ 1e-2, non-robust}
       \end{subfigure}
    \vspace{0.2cm}
    \begin{subfigure}[t]{0.3\textwidth}
        \centering
        \includegraphics[scale=0.13, trim=90pt 0 150pt 0, clip]{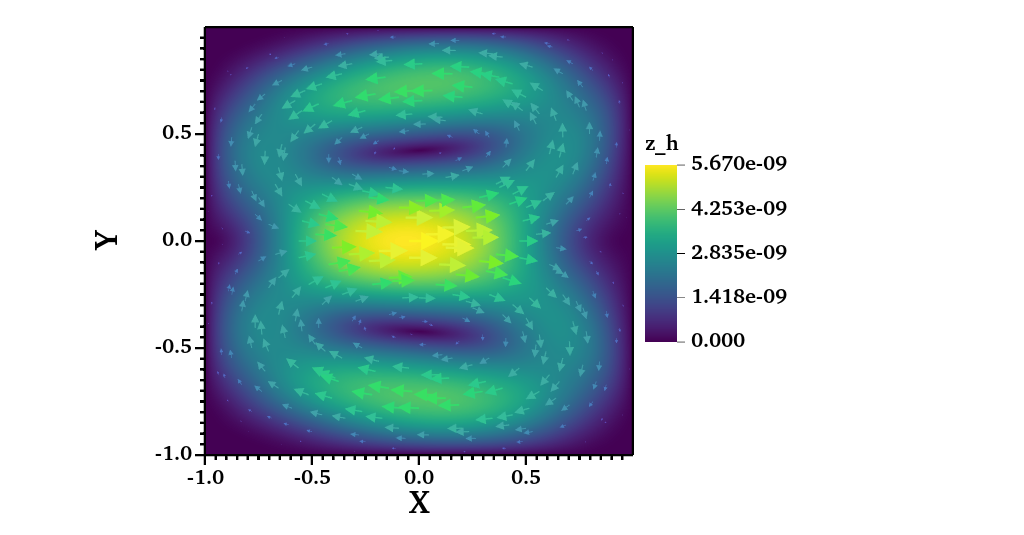}
        \caption{$\nu = $ 1e+0, robust}
      \end{subfigure}
      \hfill
      \begin{subfigure}[t]{0.3\textwidth}
        \centering
        \includegraphics[scale=0.13, trim=150pt 0 150pt 0, clip]{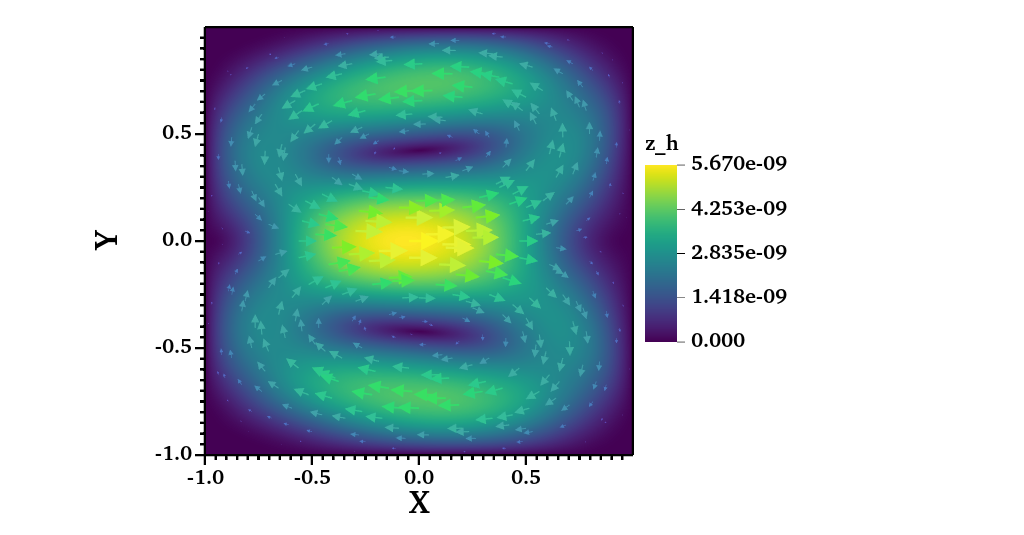}
        \caption{$\nu = $ 1e-1, robust}
      \end{subfigure}
      \hfill
      \begin{subfigure}[t]{0.3\textwidth}
        \centering
        \includegraphics[scale=0.13, trim=150pt 0 150pt 0, clip]{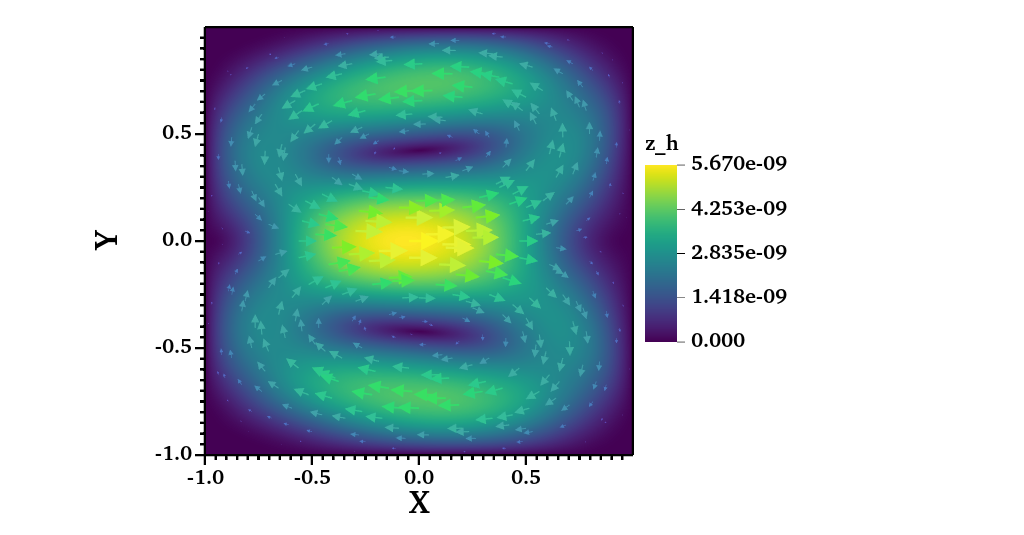}
        \caption{$\nu = $ 1e-2, robust}
      \end{subfigure}
    \caption{Magnitude and flow field of $\mathbf{z}_h$ for three different viscosities, 
    considering the non-robust (top row) and robust modification (bottom row) of the convective scheme.}
    \label{fig:adjoint}
\end{figure*}

Further, in Figure \ref{fig:adjoint} the magnitude of $\mathbf{z}_h$ along with the flow field is plotted.
The differences in accuracy between the robust and non-robust schemes highlight the effectiveness
of the interpolation operator $\pidiv$ in the schemes to restore gradient-robustness. 
This outcome aligns with the expectations and confirms the applicability of gradient
robustness in the setting of optimal control.
\hspace*{-1.7cm}


\section*{Acknowledgments}
  This work is a contribution to the project (Surface Wave-Driven Energy Fluxes At The Air-Sea
  Interface) of the Collaborative Research Center TRR 181 "Energy Transfers in Atmosphere
  and Ocean" funded by the Deutsche Forschungsgemeinschaft (DFG, German Research Foundation)
  - Projektnummer 274762653."

\section*{Competing Interests}
The authors have no conflicts of interest to declare that are relevant to the content of this chapter.


\bibliographystyle{abbrv}
\bibliography{References,dealanddope}

\end{document}